\documentclass[english,mathpazo]{elsarticle}
\usepackage[T1]{fontenc}
\usepackage[latin9]{inputenc}
\usepackage{color}
\usepackage{float}
\usepackage{amsmath}
\usepackage{amssymb}
\usepackage{graphicx}

\makeatletter

\providecommand{\tabularnewline}{\\}
\newcommand{\lyxdot}{.}

\floatstyle{ruled}
\newfloat{algorithm}{tbp}{loa}
\providecommand{\algorithmname}{Algorithm}
\floatname{algorithm}{\protect\algorithmname}

\usepackage{graphics}
\usepackage{xcolor}

\def\ps@pprintTitle{%
   \let\@oddhead\@empty
   \let\@evenhead\@empty
   \let\@oddfoot\@empty
   \let\@evenfoot\@oddfoot
}

\makeatother

\usepackage{babel}
\usepackage{listings}

\begin{document}
\title{A split step Fourier/discontinuous Galerkin scheme for the Kadomtsev--Petviashvili equation\tnoteref{label1}} \tnotetext[label1]{
The computational results presented have been achieved [in part] using the Vienna Scientific Cluster (VSC).}
\author[uibk]{Lukas Einkemmer\corref{cor1}} \ead{lukas.einkemmer@uibk.ac.at}
\author[uibk]{Alexander Ostermann} \ead{alexander.ostermann@uibk.ac.at}
\address[uibk]{Department of Mathematics, University of Innsbruck, Austria}
\cortext[cor1]{Corresponding author}
\begin{abstract} In this paper we propose a method to solve the Kadomtsev--Petviashvili equation based on splitting the linear part of the equation from the nonlinear part. The linear part is treated using FFTs, while the nonlinear part is approximated using a semi-Lagrangian discontinuous Galerkin approach of arbitrary order.

We demonstrate the efficiency and accuracy of the numerical method by providing a range of numerical simulations. In particular, we find that our approach can outperform the numerical methods considered in the literature by up to a factor of five. Although we focus on the Kadomtsev--Petviashvili equation in this paper, the proposed numerical scheme can be extended to a range of related models as well.
\end{abstract}  
\maketitle

\section{Introduction}

The Kadomtsev\textendash Petviashvili (KP) equation has been proposed
in \cite{kadomtsev1970}. It models wave propagation in dispersive
and nonlinear media and is usually stated as follows
\begin{equation}
(u_{t}+6uu_{x}+\epsilon^{2}u_{xxx})_{x}+\lambda u_{yy}=0,\label{eq:KP-equation}
\end{equation}
where $\lambda$ and $\epsilon$ are parameters. Although it is most
well known to model the propagation of water waves \cite{ablowitz1979}
it has found applications in such diverse areas as ion-acoustic plasma
waves \cite{kadomtsev1970,ohno2016}, the semiclassical description
of motion for lattice bosons \cite{demler20111}, and waves in compressible
hyperelastic plates \cite{chen2006}. From a mathematical point of
view the equation can be understood as a two-dimensional generalization
of the Korteweg\textendash de Vries equation (KdV) equation and it
shares many interesting properties, such as soliton solutions and
blow-up in finite time \cite{klein2007,klein2011,klein2012,minzoni1996},
with that system. 

It has been found that solving the KP equation poses significant challenges
for numerical methods. While it is clear that explicit schemes suffer
from a very severe stability restriction on the time step size (due
to the third order derivative in space), it was found in \cite{klein2011}
that implicit and implicit-explicit (IMEX) methods often do not converge (i.e.~that the iterative linear or nonlinear solvers necessary for the efficient implementation of these schemes do not converge).
Thus, exponential integrators have been suggested and successfully
applied to the KP equation \cite{klein2011}. However, it has been
realized recently that splitting schemes can be competitive especially
in the low to medium precision regime most of interest in applications.
The reason why splitting methods have been dismissed up to that point
is that no efficient schemes were available to solve the Burgers'
type nonlinearity. In \cite{einkemmer2015} this problem was solved
by using a semi-Lagrangian method based on Lagrange interpolation.
Nevertheless, employing Lagrange interpolation is perhaps the most
significant weakness of that algorithm. It is well known that Lagrange
interpolation is quite diffusive (see, for example, \cite{sonnendrucker1999,filbet2003})
which is certainly not ideal taking the Hamiltonian structure (i.e.~dissipation
free nature) of the KP equation into account. This is in particular
true as for the exponential integrators (pseudo) spectral methods
can be used in the nonlinear part as well.

In the context of the Vlasov equation semi-Lagrangian methods have
a long tradition. In recent years, semi-Lagrangian discontinuous Galerkin
methods have been studied extensively as these methods provide a local
approximation (which facilitates parallelization) and are only slightly
diffusive (see \textcolor{black}{\cite{qiu2011,crouseilles2011,rossmanith2011,einkemmer2015,einkemmer2016,einkemmer1211,einkemmer1602,mehrenberger2013,bokanowski2016convergence}}). However, in this case only constant advection problems have to be solved. \textcolor{black}{Position dependent problems have been considered both in the context of kinetic problems \cite{crouseilles2012guiding,qiu2011} and in the context of tracer-particle transport phenomena \cite{guo2014conservative,cai2017}. Some of these techniques have also been extended to second order equations \cite{bokanowski2016semi}. A nonlinear problem with a constant advection has been considered in \cite{bokanowski2016convergence}. Let us duly note, however, that in all of the mentioned works the advection is still linear.}

In the present paper we propose a semi-Lagrangian discontinuous
Galerkin method \textcolor{black}{for a nonlinear advection problem (more specifically, Burgers' equation)} that is used within a time splitting approach to solve
the KP equation. The proposed algorithm can take arbitrarily large
steps in time with only a very modest increase in the run time. The
algorithm as well as some computational aspects are described in section
\ref{sec:numerical-method}. We then investigate the performance of
the proposed method by conducting numerical simulations for a range
of problems (section \ref{sec:Numerical-simulations}). Finally, we
discuss some directions of future work and conclude in section \ref{sec:conclusion}.

\section{Description of the numerical method\label{sec:numerical-method}}

We consider the KP equation as stated in equation (\ref{eq:KP-equation}).
In accordance with the literature we choose either $\lambda=1$ (weak
surface tension) or $\lambda=-1$ (strong surface tension) and call
the latter the KP I model and the former the KP II model.

Before a numerical scheme is applied we rewrite equation (\ref{eq:KP-equation}),
assuming periodic boundary conditions, in evolution form
\begin{equation}
u_{t}+6uu_{x}+\epsilon^{2}u_{xxx}+\lambda\partial_{x}^{-1}u_{yy}=0,\label{eq:KP-evolution}
\end{equation}
where $\partial_{x}^{-1}$ is understood as the regularized Fourier
multiplier of $-i/k_{x}$. That is, we use the Fourier multiplier
\[
\frac{-i}{k_{x}+i\lambda\delta}
\]
with $\delta=2^{-52}$ (the machine epsilon for double precision floating
point arithmetic). \textcolor{black}{In this work, as in \cite{klein2011}, only classical solutions of the KP equation will be considered.}

As a first step, we split equation (\ref{eq:KP-evolution}) into a
linear part
\begin{equation}
u_{t}=Au=-\epsilon^{2}u_{xxx}-\lambda\partial_{x}^{-1}u_{yy}\label{eq:linear-part}
\end{equation}
and a nonlinear part
\begin{equation}
u_{t}=B(u)=-6uu_{x}.\label{eq:nonlinear-part}
\end{equation}
The former is solved using fast Fourier transform (FFT) techniques
(we denote the solution for the initial value $u^{0}$ at time $\tau$
as $\mathrm{e}^{\tau A}u^{0}$). The numerical method proposed for
the latter is described in detail in the remainder of this section
(we denote the solution for the initial value $u^{0}$ at time $\tau$
as $\varphi_{\tau}^{B}(u^{0})$). 

\textcolor{black}{
    Before proceeding, let us emphasize that even if the solution to equation (\ref{eq:KP-evolution}) is sufficiently regular, shocks (i.e.~discontinuous solutions) can develop as equation (\ref{eq:nonlinear-part}) is integrated in time. Since our algorithm for Burgers' equation is not equipped to handle such discontinuous solutions, this imposes a restriction on the step size $\tau$. Note, however, that this is \textit{not} a CFL condition as the restriction only depends on properties of the solution (and not on the grid size). In any case, numerical simulations conducted in \cite{einkemmer2015} show that for the KP equation, the problem we are ultimately interested in, this is not an issue. More precisely, the step size is dictated by accuracy constraints and we are able to take comparable or even larger time steps to what has been reported in the literature (for exponential integrators or implicit numerical methods, see for example \cite{klein2011}).
}

In the present work we will exclusively
use the second order accurate Strang splitting scheme
\[
u^{n+1}=\mathrm{e}^{\frac{\tau}{2}A}\varphi_{\tau}^{B}(\mathrm{e}^{\frac{\tau}{2}A}u^{n}),
\]
where $u^{n}$ is an approximation to the exact solution $u(t^{n})$
at time $t^{n}=n\tau$. Let us note, however, that in \cite{einkemmer2015}
extensions to fourth order splitting methods have been discussed as
well.

\subsection{The semi-Lagrangian discontinuous Galerkin approach for the nonlinearity}

Our goal in this section is to compute an approximation to the nonlinear
part of the splitting procedure (given in equation (\ref{eq:nonlinear-part})).
That is, we consider Burgers' equation
\[
\partial_{t}u(t,x)+u(t,x)\partial_{x}u(t,x)=0,\qquad\qquad u(0,x)=u^{0}(x).
\]
By using the method of characteristics (see, for example \cite{zauderer2011})
we obtain 
\begin{align*}
X^{\prime}(t;t^{n+1},x) & =U(t;t^{n},u^{0}(X(t^{n};t^{n+1},x)))\\
U^{\prime}(t;t^{n},u^{0}(x)) & =0,
\end{align*}
where $X(t;s,x)$ and $U(t;s,u)$ are the characteristic curves of
the position and the solution, respectively, that satisfy $X(s;s,x)=x$
and $U(s;s,u)=u$. Note that since $x$ is known at time $t^{n+1}$
and $u^{0}$ is known at time $t^{n}$ we have to integrate $X$ backward
in time and $U$ forward in time. However, since $U$ is constant
in time we immediately obtain
\begin{align*}
U(t;t^{n},u^{0}(x)) & =u^{0}(x)
\end{align*}
which implies that

\begin{equation}
X^{\prime}(t;t^{n+1},x)=u^{0}(X(t^{n};t^{n+1},x)).\label{eq:eq-to-solve}
\end{equation}
Therefore, we have eliminated $U$ and the desired solution of Burgers'
equation can be written as
\begin{equation}
u(t,x)=u^{0}(X(t^{n};t^{n+1},x)).\label{eq:sol-Burgers}
\end{equation}
We still have to integrate $X$ (starting from $x$) backward in time
in order to obtain the solution at $x$ (i.e.~we have to solve equation
(\ref{eq:eq-to-solve})). However, since the right hand side is constant
we can rewrite this differential equation as an algebraic equation
\begin{equation}
X(t^{n};t^{n+1},x)=x-\tau u^{0}(X(t^{n};t^{n+1},x)),\label{eq:algebraic-eq}
\end{equation}
where $\tau=t^{n+1}-t^{n}$ is the time step size.
Note that obtaining a unique solution of equation (\ref{eq:algebraic-eq}) is predicated on the assumption
that the characteristics do not cross. 

The most straightforward
approach to compute an approximate solution of (\ref{eq:algebraic-eq})
is by conducting a fixed-point iteration. That is, we perform the
following iteration
\[
\alpha^{(k+1)}=x-\tau u^{0}(\alpha^{(k)}),\qquad\qquad\alpha^{(0)}=x.
\]
Until now we have left space continuous. Thus, the only approximation
made is due to the truncation of the fixed-point iteration. However,
in order to implement a viable numerical scheme on a computer, a space
discretization strategy has to be employed. As stated in the introduction
we will perform the space discretization in the context of the semi-Lagrangian
discontinuous Galerkin (sLdG) approach. 

First, we divide our computational domain $[a,b]\subset\mathbb{R}$
into a number of cells $I_{i}=[x_{i-1/2},x_{i+1/2}]$ with size $h=x_{i+1/2}-x_{i-1/2}$.
For simplicity, we only consider an equidistant grid here. However,
the method described can be extended easily to the case of varying
cell sizes. In the $i$th cell the function values at the points $x_{ij}$
are stored, where $x_{ij}=a+h(i-1/2+\xi_{j})$ and $\xi_{j}$ is the
$j$th Gauss\textendash Legendre quadrature node scaled to the interval
$[0,1]$. The index $j$ runs from $0$ to $k$, where $k$ is the
polynomial degree that is used in each cell (note that $o=k+1$ is
the formal order of this approximation). Within this discretization
an approximation of the analytic solution is expressed as follows
\begin{equation}
u^{n}(x)=\sum_{ij}u_{ij}^{n}\ell_{ij}(x),\label{eq:approx}
\end{equation}
where $\ell_{ij}$ is the $j$th Lagrange polynomial in cell $I_{i}$
based on the nodes $x_{i0},\dots,x_{ik}$ (i.e.~the Gauss\textendash Legendre
quadrature nodes scaled and shifted to the $i$th cell). Note that
this setup implies that $u_{ij}\approx u(x_{ij})$ (i.e.~the degrees
of freedom of the numerical scheme are approximations to the function
values at the Gauss\textendash Legendre quadrature nodes in each cell). 

Now, we have to determine the coefficients at the next time step $u_{ml}^{n+1}$
from the coefficients $u_{ml}^{n}$. We start by multiplying the solution
of (\ref{eq:sol-Burgers}) with $\ell_{ml}$ and integrate in order
to obtain
\[
\int_{I_{m}}u^{n+1}(x)\ell_{ml}(x)\,\mathrm{d}x=\int_{I_{m}}u^{n}(X(t^{n};t^{n+1},x))\ell_{ml}(x)\,\mathrm{d}x
\]
which yields 

\[
\frac{h\omega_{l}}{2}u_{ml}^{n+1}=\int_{I_{m}}u^{n}(X(t^{n};t^{n+1},x))\ell_{ml}(x)\,\mathrm{d}x,
\]
where $\omega_{l}$ are the quadrature weights. That is, we perform
a projection of the analytic solution to the subspace of polynomials
up to degree $k$ (in each cell). Evaluating this integral is not
as straightforward as one might assume at first. This is due to the
fact that the function $u^{n}$ is discontinuous and a direct application
of Gauss\textendash Legendre quadrature (or any other quadrature rule
for that matter) would not give the (up to machine precision) exact
result.

To remedy this situation we proceed as follows. First, we use (\ref{eq:approx})
to write
\begin{equation}
u_{ml}^{n+1}=\frac{2}{h\omega_{l}}\int_{I_{m}}u^{n}(X(t^{n};t^{n+1},x))\ell_{ml}(x)\,\mathrm{d}x=\frac{2}{h\omega_{l}}\sum_{ij}u_{ij}^{n}\int_{I_{m}}\ell_{ij}(x^{b})\ell_{ml}(x)\,\mathrm{d}x,\label{eq:update-naive}
\end{equation}
where $x^{b}=X(t^{n};t^{n+1},x)$ is used as a shorthand to denote
endpoint of the characteristic starting at $x$. Now, only a small
minority of integrals in the sum will contribute a non-zero value
to $u_{ml}^{n+1}$. To identify the indices $i$ for which this is
true, we have to determine to what extend the support of $\ell_{ij}(x^{b})$
overlaps with $I_{m}$. There are two possibilities of nonempty intersection
that we handle separately in the implementation of our numerical scheme:
\begin{itemize}
\item Full overlap ($\text{supp}\;\ell_{ij}(x^{b})\supset I_{m}$): In this
case the function $\ell_{ij}(x^{b})$ is polynomial in $I_{m}$ and
we can apply a Gauss\textendash Legendre quadrature rule in order
to compute the integral exactly;
\item Partial overlap: In this case the function $\ell_{ij}(x^{b})$ has
one or two discontinuities in the interval $I_{m}$. In this case
we apply a Gauss\textendash Legendre quadrature to $\text{supp}\;\ell_{ij}(x^{b})\cap I_{m}$.
\end{itemize}
To determine which of these cases apply (and in the latter case what
the endpoints of integration should be) we follow the characteristics
forward in time. That is, to obtain $\text{supp}\;\ell_{ij}(x^{b})=[a_{i},b_{i}]$
we have to compute
\[
a_{i}=X(t^{n+1};t^{n},x_{i-1/2}),\qquad\qquad b_{i}=X(t^{n+1};t^{n},x_{i+1/2}).
\]
Following the characteristics forward in time is less involved since
we now have to solve
\begin{align*}
X'(t;t^{n},x) & =u^{0}(X(t^{n};t^{n},x))
\end{align*}
which immediately yields
\[
X(t^{n+1};t^{n},x)=x+\tau u^{n}(x).
\]
Thus, we have 
\begin{equation}
a_{i}=x_{i-1/2}+\tau u^{n}(x_{i-1/2}),\qquad\qquad b_{i}=x_{i+1/2}+\tau u^{n}(x_{i+1/2}).\label{eq:predict-1}
\end{equation}
The difficulty here is that $u^{n}(x_{i-1/2})$ and $u^{n}(x_{i+1/2})$
are not well defined due to the discontinuous approximation used for
$u^{n}(x)$. Thus, we replace equation (\ref{eq:predict-1}) by
\begin{equation}
a_{i}=x_{i-1/2}+\tau\tilde{u}^{n}(x_{i-1/2}),\qquad\qquad b_{i}=x_{i+1/2}+\tau\tilde{u}^{n}(x_{i+1/2})\label{eq:predict}
\end{equation}

with
\[
\tilde{u}_{i-1/2}^{n}=\frac{u^{n}(x_{i-1/2}+)+u^{n}(x_{i-1/2}-)}{2},
\]
here we have used $+$ and $-$ to denote the right-sided and left-sided
limit, respectively. The choice made yields the shock speed consistent
with the exact solution of the Riemann problem for Burgers' equation.
Note that although the discontinuity is on the order of the approximation
error and thus the accuracy of the numerical scheme is not negatively
impacted, we will see in section \ref{subsec:conservation-mass} that
conservation of mass is lost.

\begin{figure}
\begin{centering}
\includegraphics[width=8cm]{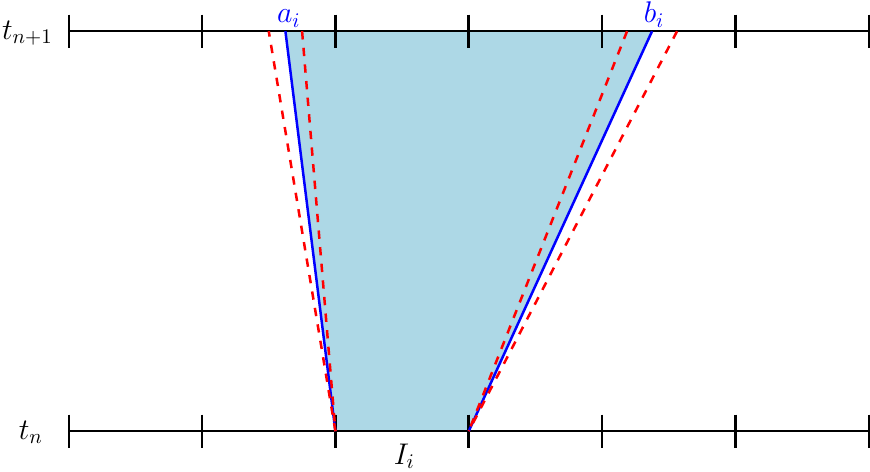}
\par\end{centering}
\begin{centering}
\includegraphics[width=8cm]{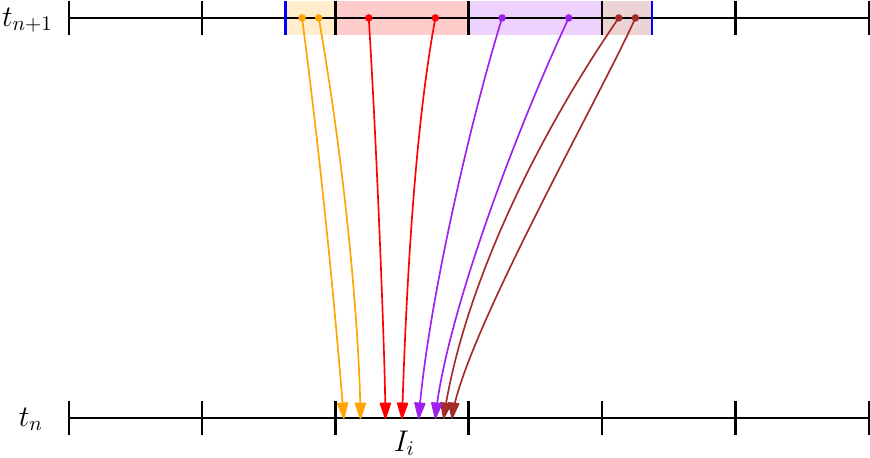}
\par\end{centering}
\caption{An illustration of the proposed algorithm is shown. First, a forward
step (shown on the top) is conducted that determines the domain of
influence of the cell $I_{i}$ (the blue lines). Note that as has
been pointed out in the text, the true characteristics (dashed red
lines) might cross and thus this only yields an approximation. Then
in the backward step (shown on the bottom) the region between the
two blue lines is divided into parts (illustrated with different colors)
each of which overlap with a single cell. In these parts we integrate
the characteristics back in time starting at the corresponding Gauss\textendash Legendre
quadrature points. The obtained values are then used to update the
corresponding degrees of freedom. This backward step is exact up to
the error committed in integrating the characteristics. \label{fig:illustration-algorithm}}
\end{figure}

In principle, we could use equation (\ref{eq:update-naive}) as the
basis of our numerical implementation. However, it is beneficial to
consider
\begin{equation}
u_{ml}^{n+1}=\frac{2}{h\omega_{l}}\sum_{i}\int_{I_{m}}u_{i}^{n}(x^{b})\ell_{ml}(x)\,\mathrm{d}x=\frac{2}{h\omega_{l}}\sum_{i}\int_{I_{m}\cap I_{i}}u_{i}^{n}(x^{b})\ell_{ml}(x)\,\mathrm{d}x\label{eq:update}
\end{equation}
where $u_{i}^{n}(x)=\sum_{j}u_{ij}^{n}\ell_{ij}(x)$. That is, $u_{i}^{n}(x)$
is equal to $u^{n}(x)$ in the cell $I_{i}$ and zero outside. Let
us assume that $I_{m}\cap I_{i}=[c,d]$ then we have
\begin{equation}
\int_{[c,d]}u_{i}^{n}(x^{b})\ell_{ml}(x)\,\mathrm{d}x\approx\frac{d-c}{2}\sum_{r}\omega_{r}u_{i}^{n}((\xi_{r})^{b})\ell_{ml}(\xi_{r}),\label{eq:update-discr}
\end{equation}
where $\xi_{r}$ are the Gauss\textendash Legendre quadrature points
in the interval $[c,d]$. In the actual implementation we perform
the algorithm in reverse to what (\ref{eq:update}) would suggest.
That is, instead of setting one $u_{ml}^{n+1}$ at a time (using multiple
$u_{i}^{n}$ in the process) we start with a fixed cell $i$ and then
use $u_{i}^{n}$ to update multiple cells (i.e.~update $u_{ml}^{n+1}$
for multiple $m$). This is illustrated with the pseudocode in Algorithm
\ref{alg:pseudocode} and the pictures in Figure \ref{fig:illustration-algorithm}.

\begin{algorithm}
\begin{lstlisting}[tabsize=4]
loop over all cells i
	determine a_i,b_i from (2.6)
	loop over all cells m that overlap with [a_i,b_i]
		determine the interval of quadrature [c,d]
		loop over the Gauss-Legendre points in [c,d]
			follow characteristics backward in time
		loop over all nodes l in the cell m
			update u_{ml}^{n+1} using (2.7) and (2.8)
\end{lstlisting}

\caption{Pseudocode that illustrates the implementation of the semi-Lagrangian
discontinuous Galerkin method for Burgers' equation. \label{alg:pseudocode}}
\end{algorithm}

\subsection{Computational complexity}

The pseudocode in Algorithm \ref{alg:pseudocode} gives us a convenient
vehicle to analyze the computational complexity of the present algorithm.
In the second loop we iterate over all cells that have a nonempty
overlap with the interval $[a_{i},b_{i}]$. In the constant advection
case (see, for example, \cite{qiu2011,crouseilles2011,rossmanith2011,einkemmer2016,einkemmer1211,einkemmer1602,mehrenberger2013,einkemmer2015hpc})
this implies that we have to consider at most two adjacent cells.
However, in the present case the number of cells we have to consider
depends on $\tilde{u}^{n}(x_{i+1/2})$ and $\tilde{u}^{n}(x_{i-1/2})$.
That is, it depends on how much the magnitude of the solution changes
from one cell interface to the next. We note, however, that even for
large CFL numbers it is uncommon that more than three cells are processed
in that step. Since, except for the first loop (which iterates over
the number of cells $n$), the remaining loops only iterate over nodes
within a single cell, the numerical algorithms scales as $\mathcal{O}(n)$,
i.e.~linearly in the number of cells $n$.

Let us also consider what happens as we increase the order $o$ of
the numerical scheme. The most computationally intensive part of the
algorithm is, by a significant margin, to follow the characteristics
backward in time. This function is called $o$ times (once for each
quadrature point) and requires a single function evaluation of $u^{n}$
for each iteration of the fixed-point algorithm. Each of these function
evaluations scales linearly in the order of the method. Thus, in total
the numerical algorithms scales as $\mathcal{O}(o^{2})$. Thus, strictly
speaking we do not have linear scaling in the degrees of freedom.
However, let us note that on modern computer systems the algorithm
(except for extremely large orders) is certainly memory bound. This
is further reason for choosing the implementation illustrated in Figure
\ref{alg:pseudocode} as it facilitates local memory access (both
for reading and writing data).

\subsection{Improving the fixed-point iteration\label{subsec:improv-fp}}

As mentioned in the previous section the most compute intensive part
of the algorithm is to follow the characteristics backward in time.
Thus, decreasing the number of iterations required for the fixed-point
iteration would accordingly decrease the time required to run our
algorithm. Our goal is still to solve the following fixed-point problem
\[
\alpha=x-\tau u^{0}(\alpha).
\]
Rewriting this as a root finding problem we obtain

\[
f(\alpha)=x-\tau u^{0}(\alpha)-\alpha=0
\]
to which we can apply Newton's method
\[
\alpha^{(k+1)}=\alpha^{(k)}+\frac{x-\tau u^{0}(\alpha^{(k)})-\alpha^{(k)}}{1+\tau(u^{0})^{\prime}(\alpha^{(k)})},\qquad\alpha^{(0)}=x.
\]
One drawback of this method is that we have to compute $(u^{0})^{\prime}$.
Although $(u^{0})^{\prime}$ can be precomputed and consequently the
cost is independent of the number of iterations conducted, this is
still expensive if only a few iterations are required. In addition,
it increases the code complexity significantly. However, the main
drawback is that we now require two function evaluations ($u^{0}(\alpha^{(k)})$
and $(u^{0})^{\prime}(\alpha^{(k)})$) per iteration. Nevertheless
as we will see in the next section especially for large time steps
Newton's method is clearly superior to performing a fixed-point iteration.

We can, however, remove all of the deficiencies of Newton's method,
while still retaining the speed of convergence, by employing the secant
method. In this case we start with a fixed-point iteration
\[
\alpha^{(1)}=x-\tau u^{0}(x),\qquad\alpha^{(0)}=x
\]
and then continue using the secant method
\begin{align*}
g^{(k)} & =x-\tau u^{0}(\alpha^{(k)})-\alpha^{(k)},\\
\alpha^{(k+1)} & =\alpha^{(k)}-\frac{\left(\alpha^{(k)}-\alpha^{(k-1)}\right)g^{(k)}}{g^{(k)}-g^{(k-1)}}.
\end{align*}
Within this approach we can reuse the function evaluation from the
previous iteration. Moreover, the speed of convergence compared to
the fixed-point iteration is greatly improved. We will benchmark these
three approaches in section~\ref{sec:Numerical-simulations}. 

\subsection{Conservation of mass\label{subsec:conservation-mass}}

Both the linear part (\ref{eq:linear-part}) and the nonlinear part
(\ref{eq:nonlinear-part}) conserve mass. For the latter this can
be seen most easily by writing it as a conservation law
\[
\partial_{t}u(t,x)+\partial_{x}(\tfrac{1}{2}u^{2}(t,x))=0.
\]
This is still true if we replace the exact solution of the linear
part by an approximation based on Fourier techniques. Now, let us
briefly analyze conservation of mass for the nonlinear part. In previous
work for advections that are independent of $u$ (for the case of
position dependent advection see \cite{crouseilles2012guiding}; for
the case of constant advection see, for example, \cite{qiu2011,crouseilles2011,rossmanith2011,einkemmer2015,einkemmer2016,einkemmer1211,einkemmer1602,mehrenberger2013})
mass is preserved by the semi-Lagrangian discontinuous Galerkin scheme.
However, this is not true in the present case as the approximation
of $u$ is discontinuous at the cell interface. Although that jump
is on the order of the discretization error (and thus does not negatively
impact accuracy) it implies that either the characteristics cross
(resulting in a shock wave) or fan out (resulting in a rarefaction
wave). Note that in both cases the solution between the two red dashed
lines in Figure \ref{fig:illustration-algorithm} can not be (solely)
determined by the method of characteristics. Furthermore, the solution
is, in general, not smooth across this interface. Since this, in particular,
implies that the quadrature is not exact in this region, the numerical
scheme described is not mass conservative up to machine precision.
It was already pointed out, in a slightly different context, that
continuity across the cell interface for the advection speed is necessary
in order to obtain a mass conservative semi-Lagrangian discontinuous
Galerkin scheme \cite{crouseilles2012guiding}.

In principle, this deficiency can be remedied by introducing an additional
interval (between the two red dashed lines in Figure \ref{fig:illustration-algorithm})
on which we use the analytic solution of the Riemann problem to update
the numerical solution. This adds additional computational cost but
might still be beneficial for problems where, for example, long time
integration is essential. We consider this as future work.

Strictly speaking there is one additional error made with respect
to mass conservation. Since our scheme is based on the method of characteristics
and we therefore do not approximate Burgers' equation in its conservative
form, conservation of mass is only guaranteed if the characteristics
are solved exactly. However, in our opinion, this is only a minor
issue as the numerical simulations in section \ref{sec:Numerical-simulations}
demonstrate that using the secant or Newton's method only a few iterations
are sufficient to obtain the characteristics up to machine precision. 

\section{Numerical simulations\label{sec:Numerical-simulations}}

The goal of this section is to both validate the implementation as
well as to demonstrate its efficiency for the KP equation. We start
in section \ref{subsec:numerical-burgers} by considering only the
nonlinear part of the splitting approach (i.e.~Burgers' equation).
Then in section \ref{subsec:numerical-KP} we will consider a number
of different problems in the context of the KP equation. Finally,
in section \ref{subsec:comparison} we perform a comparison of the
proposed algorithm with the exponential type methods that are considered
in \cite{klein2011}.

\subsection{Burgers' equation\label{subsec:numerical-burgers}}

In this section we evaluate the accuracy and efficiency of the proposed
algorithm for the Burgers' equation by comparing it to a spectral
implementation. More specifically, in the latter case we use the fast
Fourier transform (FFT) to perform the space discretization and solve
the resulting ordinary differential equation using the CVODE library
where both the absolute and relative tolerance are set to $10^{-15}$.
Of course, if we are only interested in solving Burgers' equation
this would mean that the spectral implementation (which performs many
substeps as CVODE implements the fully implicit BDF methods) is orders
of magnitude slower compared to our implementation which (as we will
see in this section) requires at most a moderate number of iterations.
Thus, in this context this is certainly not a sensible comparison.
However, for the application we are interested in and which we will
consider in the next section, i.e.~solving the Kadomtsev\textendash Petviashvili
equation, this is still a good indicator of the space discretization
error.

Let us start by considering Burgers' equation on the interval $[0,2\pi]$
using periodic boundary conditions and the initial value
\[
u^{0}(x)=\sin x.
\]
For the remainder of this section we will refer to this as the sine
initial value. The advantage of this choice is that the analytic solution
can be expresses as follows
\[
u(t,x)=-\sum_{k=0}^{\infty}\frac{2J_{k}(-kt)}{kt}\sin(kx),
\]
where $J_{k}(x)$ are the Bessel functions of the first kind. This
analytic solution is used to verify our implementation. Let us further
remark that the initial value $u(0,x)$ is perfectly smooth but that
at $t=1$ the solution develops a discontinuity. Thus, as we increase
$t$ we obtain a progressively less regular solution which is more
challenging for the higher order methods (including the spectral method).

\begin{figure}[H]
\begin{centering}
\includegraphics[width=6cm]{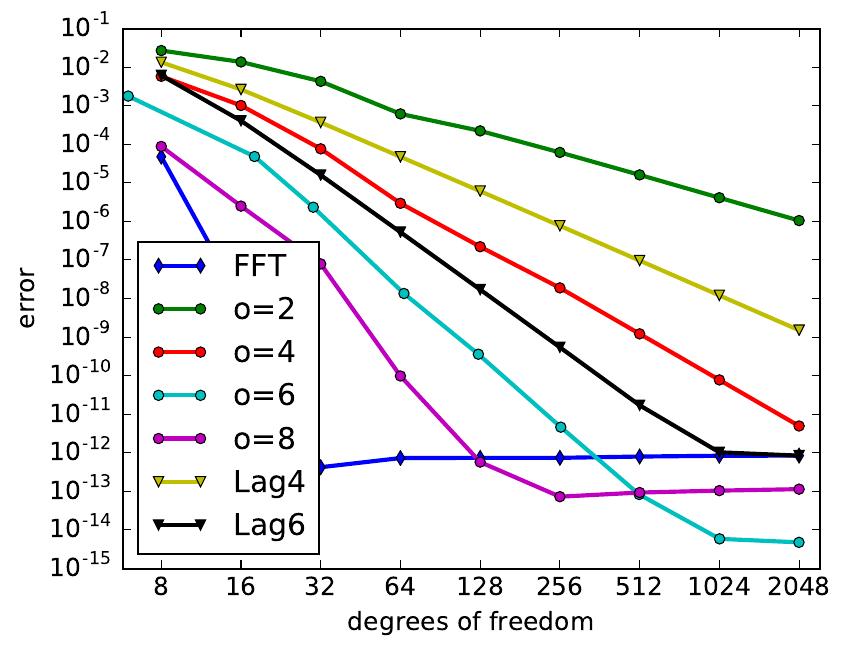}\includegraphics[width=6cm]{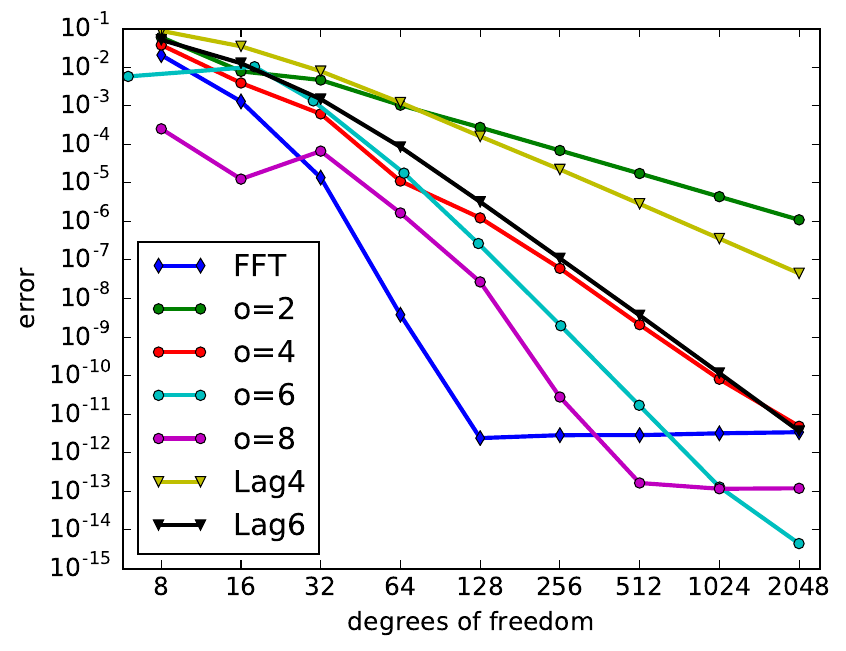}
\par\end{centering}
\begin{centering}
\includegraphics[width=6cm]{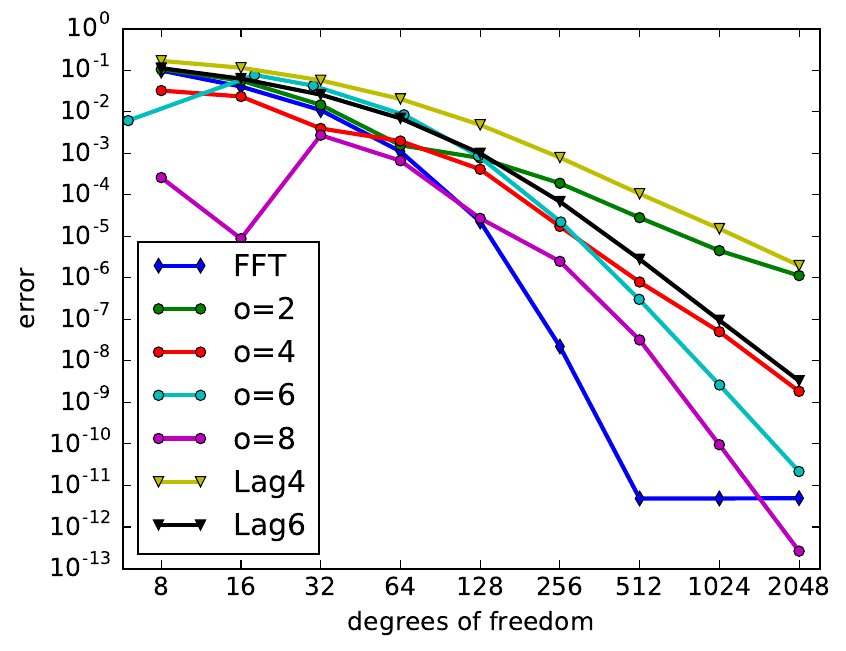}\includegraphics[width=6cm]{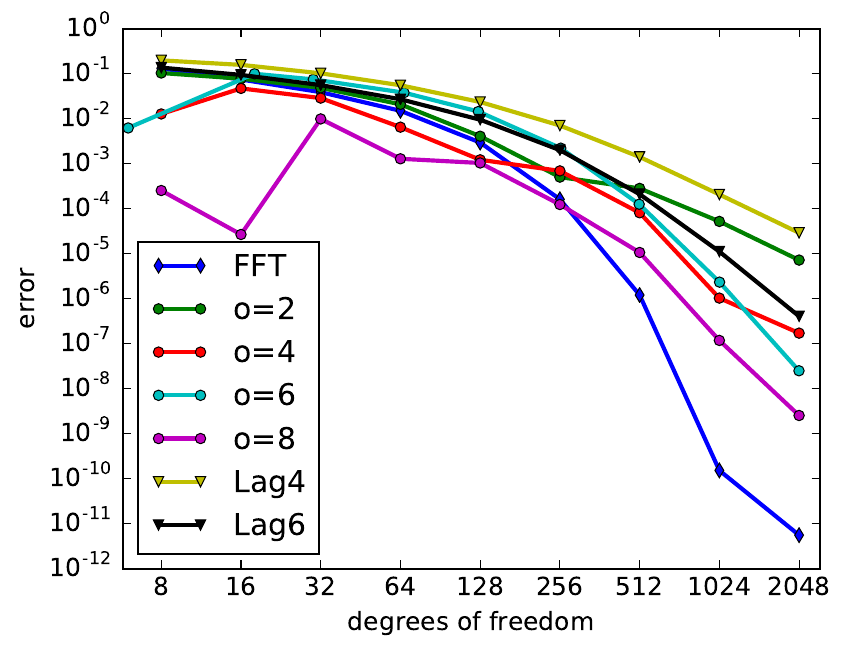}
\par\end{centering}
    \caption{\color{black} The error in the infinity norm as a function of the degrees of freedom
is shown for the sine initial value. From top left to bottom right
the error is evaluated at $t=0.1$, $0.5$, $0.8$, $0.9$ corresponding
to progressively less regular solutions. For all simulations $10$
iterations with the secant method are conducted and the numerical solution is compared to the analytic solution in order to determine the error. \label{fig:order-sine}}
\end{figure}

The numerical results for four different times are shown in Figure
\ref{fig:order-sine}. We remark that the FFT based approach is generally
superior for very regular solutions. However, in all other cases,
the eight order discontinuous Galerkin method is on par or superior
for the tolerances that are usually of interest in practical applications
(i.e.~an error at or above $10^{-5}$). We also remark that the less
regular the solution becomes the smaller the advantage of the eight
order method is compared to the fourth order method. Nevertheless,
some advantage with respect to accuracy persists even for $t=0.9$.
We have also tried using discontinuous Galerkin methods of order higher
than eight. However, those do not result in any additional improvement
and we will therefore not consider them further.

\textcolor{black}{In addition, we compare the discontinuous Galerkin method proposed in this paper to Lagrange interpolation of order 4 (as proposed in \cite{einkemmer2015}) and order 6. These results are also shown in Figure \ref{fig:order-sine}. We observe that in all configurations considered here the discontinuous Galerkin method outperforms Lagrange interpolation, of the same order, by a significant margin.}

As a second example, we will consider Burgers' equation on the interval
$[-L,L]$ with periodic boundary conditions and the initial value
\[
u^{0}(x)=\text{sech}\,x=\frac{1}{\text{cosh}\,z}.
\]
For the remainder of this section we will refer to this as the sech
initial value. Choosing this initial value is motivated by the fact
that in many applications of the KP equation we consider initial values
that decay to zero as $x\to\infty$ but are not exactly periodic on
the truncated domain that we have to employ in a numerical simulation.
This then usually means that we have to choose $L$ relatively large
in order to render errors due to the truncation of the physical domain
negligible.

\begin{figure}[H]
\begin{centering}
\includegraphics[width=6cm]{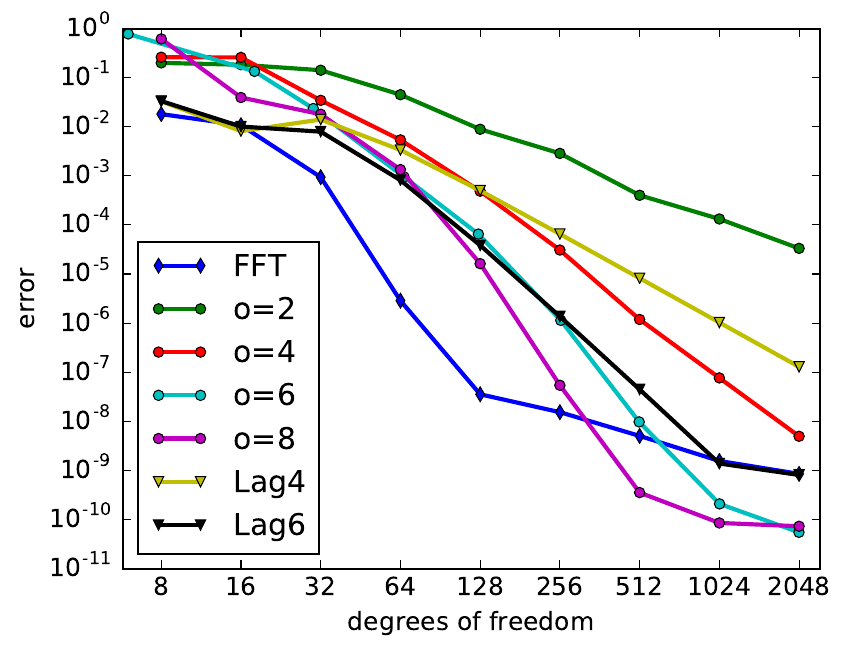}\includegraphics[width=6cm]{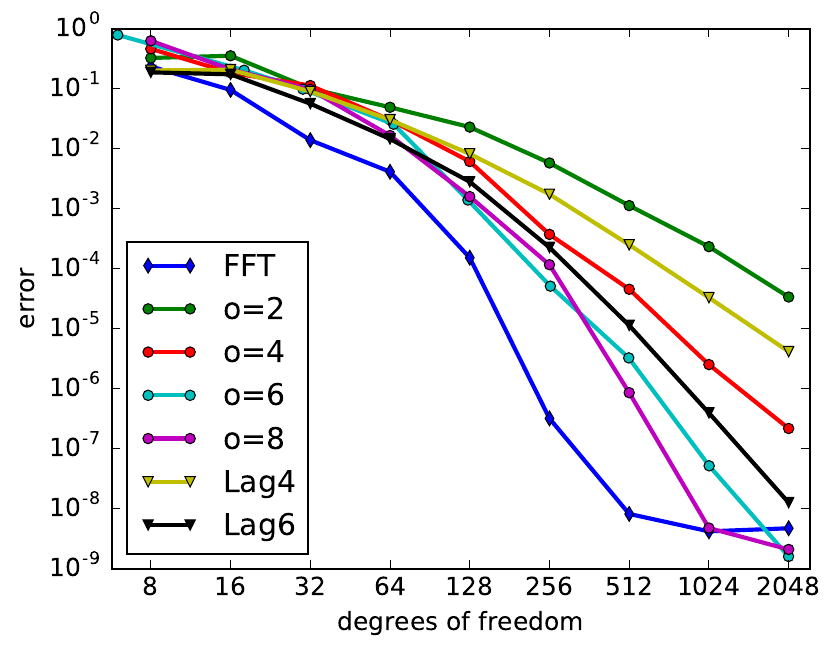}
\par\end{centering}
    \caption{\color{black} The error in the infinity norm as a function of the degrees of freedom is shown for the sech initial value. We consider a final time $t=0.2$ (left) and $t=1.1$ (right). In both cases the domain is truncated at $L=20$. For all simulations $10$ iterations with the secant method are conducted. The error is computed by comparing the numerical solution with a FFT based reference solution that uses $8192$ degrees of freedom.
\label{fig:order-sech}}
\end{figure}

\textcolor{black}{
In Figure \ref{fig:order-sech} we investigate the performance of
the discontinuous Galerkin approach as a function of regularity
(i.e.~as a function of the final time). We observe that the FFT implementation
is, in general, superior to the discontinuous Galerkin approach. However, the difference in the number of degrees of freedom needed, for a given accuracy, for the 8th order discontinuous Galerkin method is always less than a factor of two. The discontinuous Galerkin approach outperforms Lagrange interpolation (except at tolerances below $10^{-2}$ and $t=0.2$).
}

\begin{figure}[H]
\begin{centering}
\includegraphics[width=6cm]{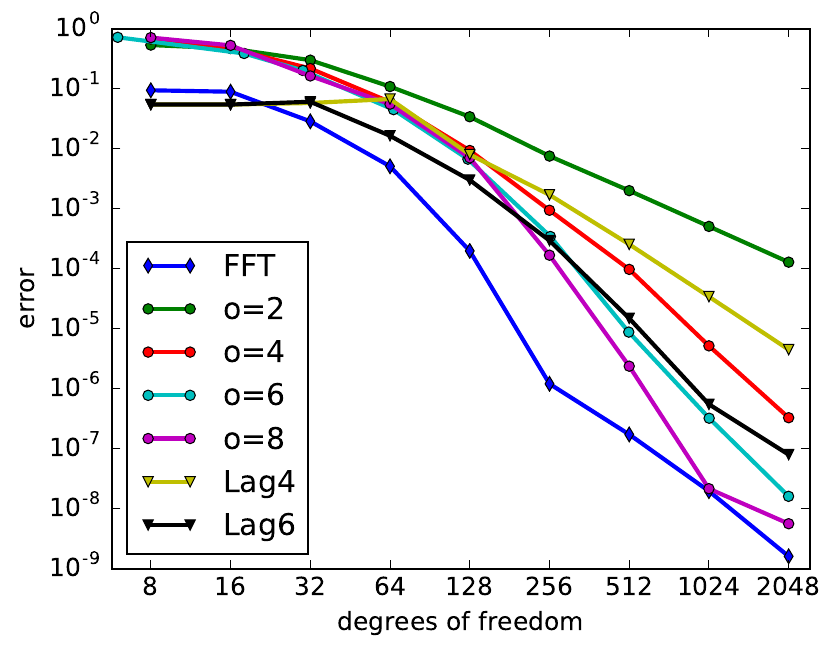}\includegraphics[width=6cm]{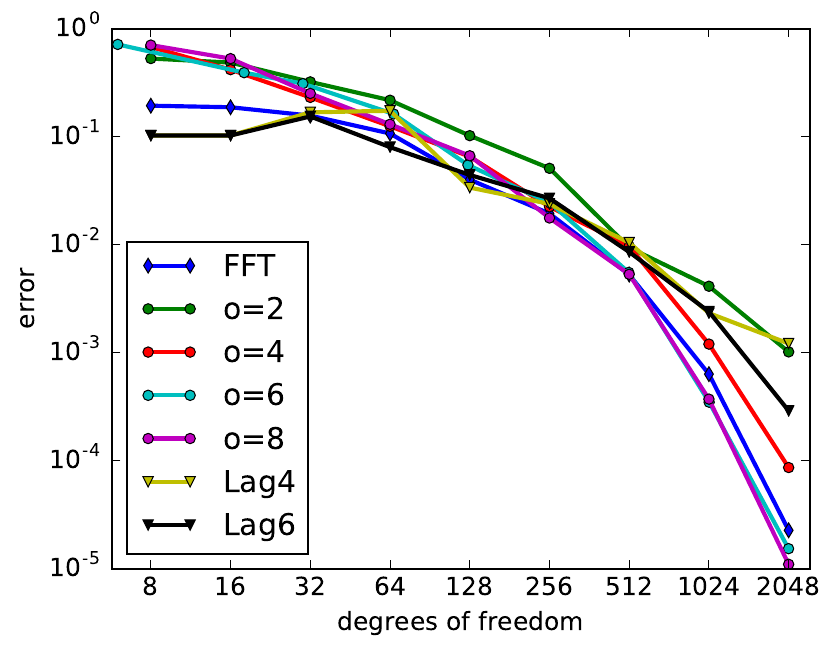}
\par\end{centering}
    \caption{\color{black} The error at $t=0.05$ (left) and $t=0.1$ (right) as a function of the degrees of freedom is shown for the oscillating initial value. For all simulations $10$ iterations with the secant method are conducted. The error is computed by comparing the numerical solution with a FFT based reference solution that uses $8192$ degrees of freedom.
 \label{fig:order-oscillating}}
\end{figure}

\textcolor{black}{
As the third and final example we consider Burgers' equation on the
interval $[-1,1]$ with periodic boundary conditions and initial value
\[
    u^{0}(x)=\cos\left(\omega (\tfrac{1}{2}+\vert x \vert) x\right) \Psi(x)
\]
with the bump function $\Psi(x)=\exp(-1/(1-x^2))$ and $\omega=20$.
For the remainder of this section we will refer to this as the oscillating
initial value. 
Once again we investigate the accuracy
of the discontinuous Galerkin and FFT implementation as a function
of the regularity of the solution. The corresponding numerical results
are shown in Figure \ref{fig:order-oscillating}. As the regularity decreases the relative performance of the eight order discontinuous Galerkin scheme, compared to the FFT implementation, increases. For the configuration with $t=0.1$ the performance of these two methods is almost identical. As before, we observe that the discontinuous Galerkin method outperforms Lagrange interpolation by a significant margin.}

The issue of domain truncation (as considered for the sech initial value) is particularly important for the KP
(and other dispersive) equation and it thus deserves some additional
comments. While we have shown here that the
discontinuous Galerkin method can be very competitive compared to the FFT
based implementation, this is certainly not the end of the story for
this particular method. What can be done relatively easily for the
discontinuous Galerkin method is to employ a grid where not all cells
have the same size. Therefore, in the spirit of a block-structured
mesh refinement scheme we would keep the present resolution in the
center of the domain (where the majority of the dynamics happens)
and use larger cells (with possible lower order) in the part of the
domain where the solution is small. This would significantly reduce
the computational effort with almost no loss in accuracy. Such an
optimization is very difficult to do within the FFT framework and
would thus further work in favor of the numerical scheme proposed
in this paper. We consider this as future work.

\begin{figure}
\begin{centering}
\includegraphics[width=6cm]{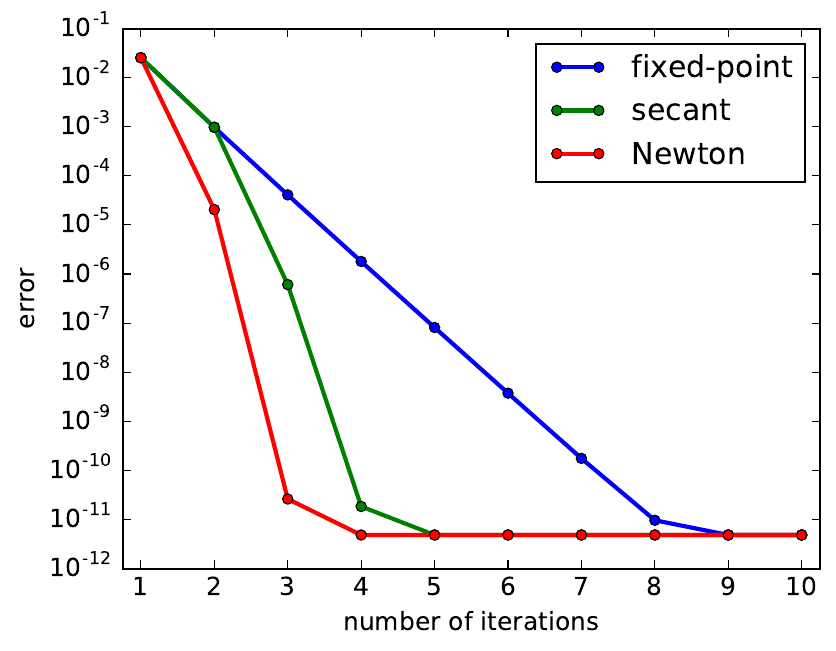}\includegraphics[width=6cm]{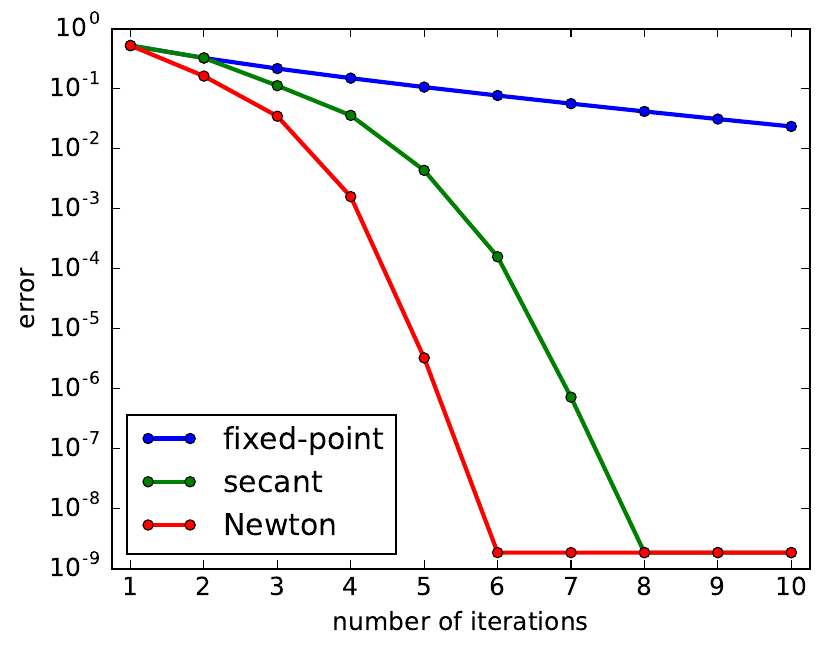} 
\par\end{centering}
\caption{Error in the infinity norm as a function of the number of iterations
is shown for the sine initial value at $t=0.05$ (left) and $t=0.8$
(right). The space discretization uses the fourth order method ($o=4$)
    with $512$ cells \textcolor{black}{and the numerical solution is compared to the analytic solution in order to obtain the error}. \label{fig:iteration-sine}}
\end{figure}

To conclude this section, let us investigate how many iterations are
required in order to obtain an accurate solution. For that purpose
we consider the three methods described in section \ref{subsec:improv-fp}
(fixed-point iteration, secant method, and Newton's method) and solve
Burgers' equation using the sine initial value. The numerical results
are shown in Figure~\ref{fig:iteration-sine}. We observe that, while
the fixed-point iteration can converge slowly if we take excessively
large time steps, for the secant method and Newton's method no more
than a couple of iterations are necessary. In fact, for these methods
the number of iterations required to reach a certain precision is
only weakly dependent on the time step size. Let us note that the
secant method requires the same computational effort compared to fixed-point
iteration while Newton's method (which requires the construction of
the derivative) is approximately twice as costly. Thus, in order to
minimize the computational effort for a given accuracy we have chosen
to exclusively use the secant method in the examples that are presented
in the next section.

\subsection{Kadomtsev\textendash Petviashvili equation\label{subsec:numerical-KP}}

First let us consider numerical results for the KP I and KP II equations
using the Schwartzian initial value
\[
u(0,x,y)=-\partial_{x}\text{sech}^{2}\left(\sqrt{x^{2}+y^{2}}\right).
\]
In all simulations we have chosen $\epsilon=0.1$ and use the computational
domain $[-10\pi,10\pi]\times[-10\pi,10\pi]$. This is precisely the
initial value that has been used in \cite{klein2011}. If not indicated
otherwise we always use $5$ iterations with the secant method to
determine the characteristic curves required in our algorithm.

\begin{figure}
\begin{centering}
\includegraphics[width=12cm]{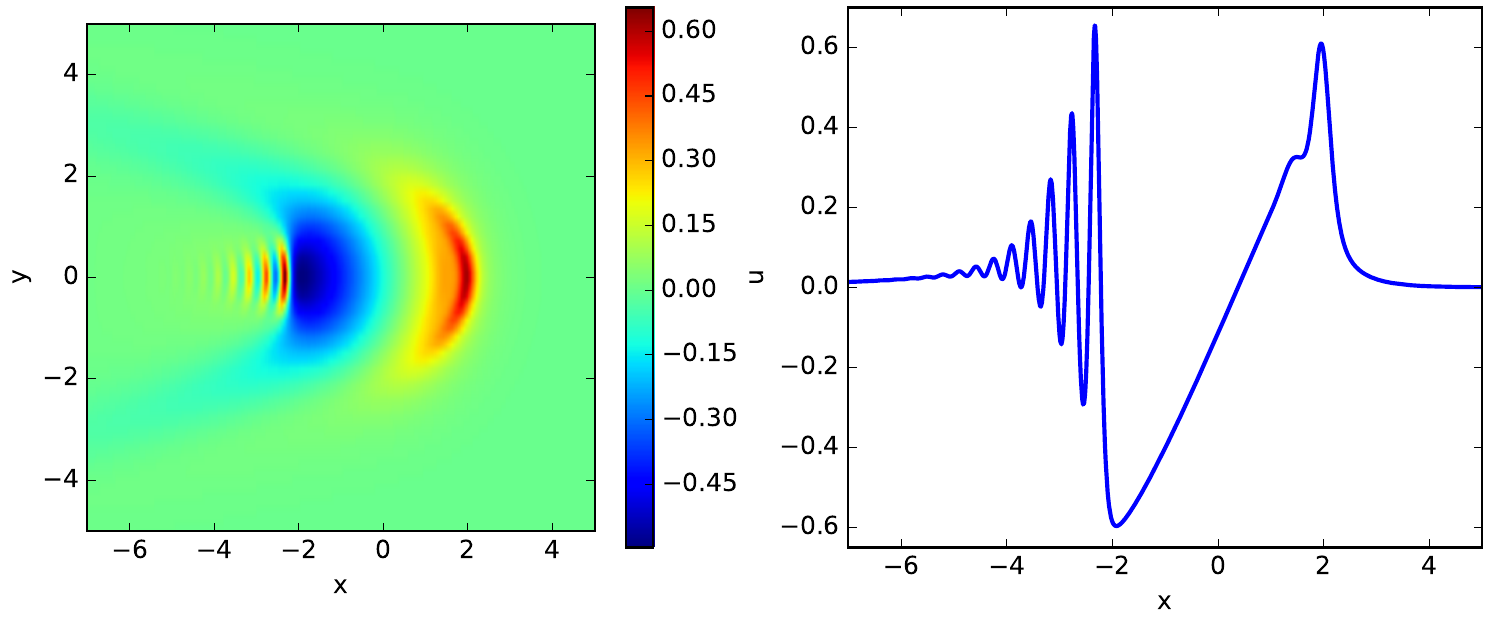}
\par\end{centering}
\begin{centering}
\includegraphics[width=6cm]{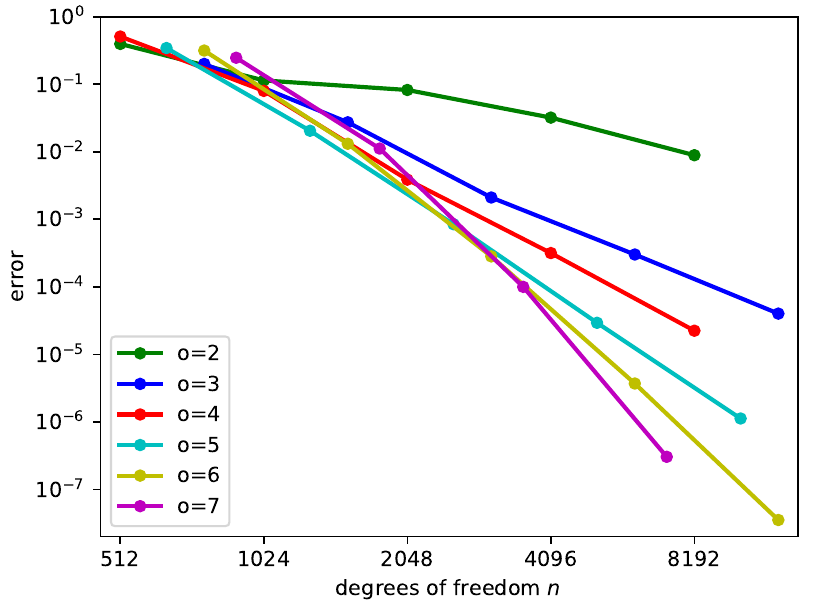}\includegraphics[width=6cm]{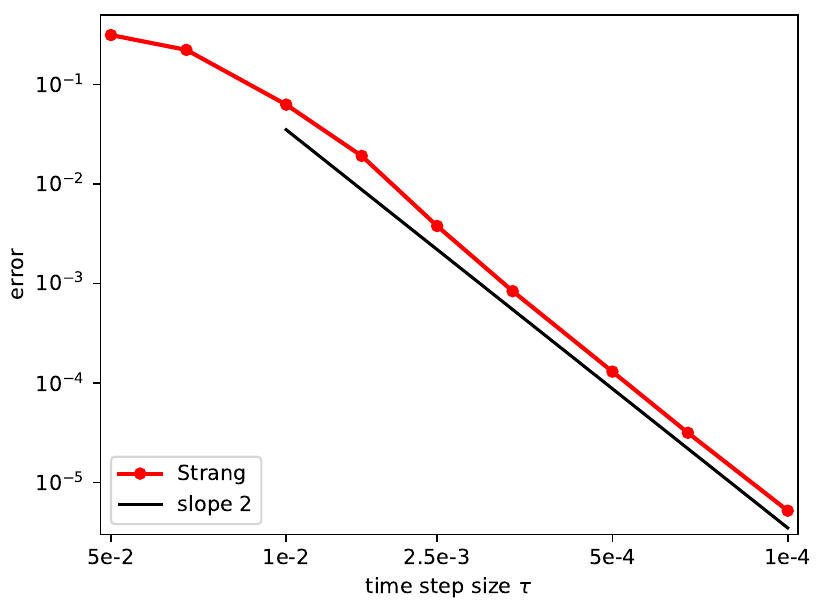}
\par\end{centering}
\caption{The numerical simulation for the KP II equation using the Schwartzian
initial value is shown at $t=0.4$ for the two-dimensional domain
(top-left) and for a one-dimensional slice with $y=0$ (top-right).
The fourth order method (for Burgers' equation) with a time step size
of $\tau=10^{-3}$ has been used. In addition, the error for the discontinuous Galerkin methods is shown
as a function of the degrees of freedom (bottom-left) and the error
    as a function of the time step size $\tau$ is shown (bottom-right). \textcolor{black}{For the former the error is determined by comparing the numerical approximation to a reference solution with $o=7$ and $28672$ degrees of freedom in the $x$-direction, while for the latter a reference solution with time step size $\tau=5 \cdot 10^{-5}$ is used.}
 \label{fig:KPII-sol-order}}
\end{figure}

First, let us consider the KP II equation (i.e.~$\lambda=1$). The
numerical results for final time $t=0.4$ are shown in Figure \ref{fig:KPII-sol-order}.
This solution agrees very well with what has been reported in the
literature (see, for example, \cite{klein2012}). In addition, we
investigate the error as the number of degrees of freedom is increased
as well as the error as we decrease the time step size (all errors
are computed in the infinity norm and are compared to a reference
solution with a sufficiently fine resolution). We find that for low
to medium precision requirements the fifth and seventh order discontinuous
Galerkin method yields the best performance, respectively. 

Before proceeding, let us discuss the number of iterations $i$ required
in the nonlinear solver. For the secant method it is generally sufficient
to perform $i=3$ iterations. The same level of precision is achieved
for the fixed-point iteration for $i=5$ iterations. We note that
the advantage of the secant method increases significantly for more
stringent tolerances. Thus, we conclude that while the relative difference
is not as large as in the numerical simulations conducted in section
\ref{subsec:numerical-burgers} (due to the fact that the splitting
error limits the time step size for a given accuracy), we can still
obtain a significant saving in computational resources by employing
the secant method.

\begin{figure}
\begin{centering}
\includegraphics[width=6cm]{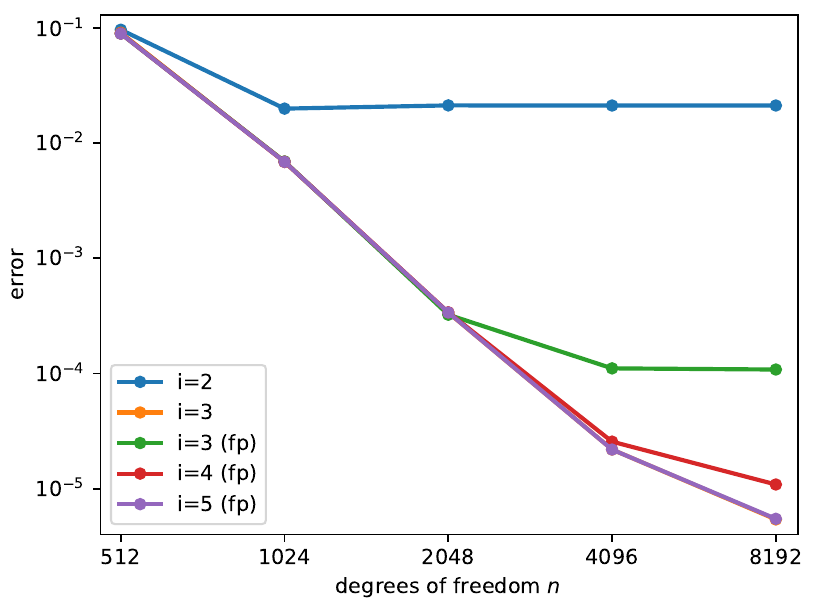}
\par\end{centering}
\caption{The error is shown as a function of the degrees of freedom for various
numbers of iterations $i$ for the fixed-point (fp) and the secant
method (note that the lines for the secant method with $i=3$ and
the fixed-point iteration with $i=5$ are indistinguishable in the
    plot). \textcolor{black}{The error is determined by comparing the numerical approximation to a reference solution with $i=20$ iteration using the secant method.} \label{fig:KPII-cost-iter}}
\end{figure}

\begin{figure}
\begin{centering}
\includegraphics[width=12cm]{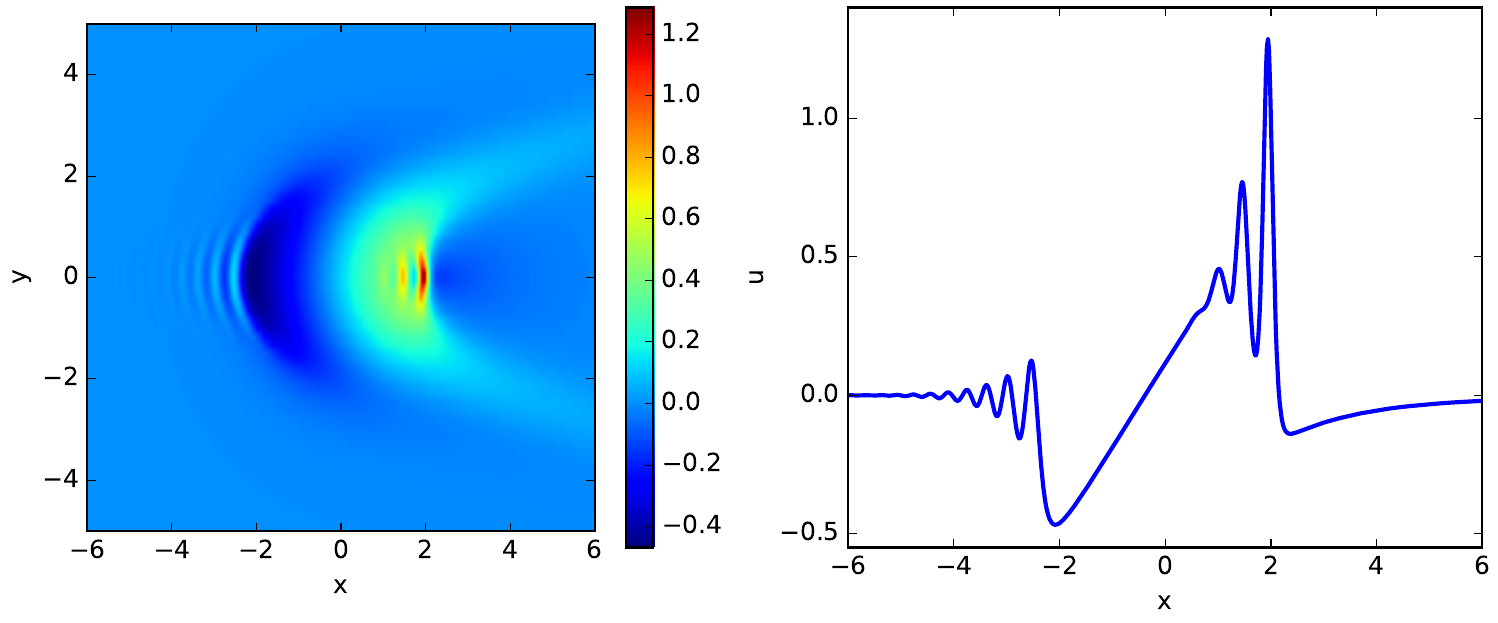}
\par\end{centering}
\centering{}\includegraphics[width=6cm]{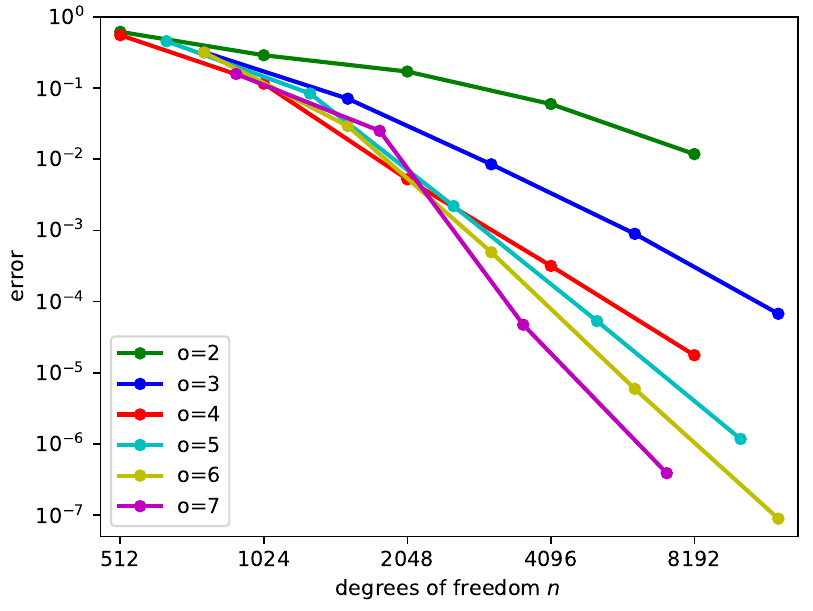}\includegraphics[width=6cm]{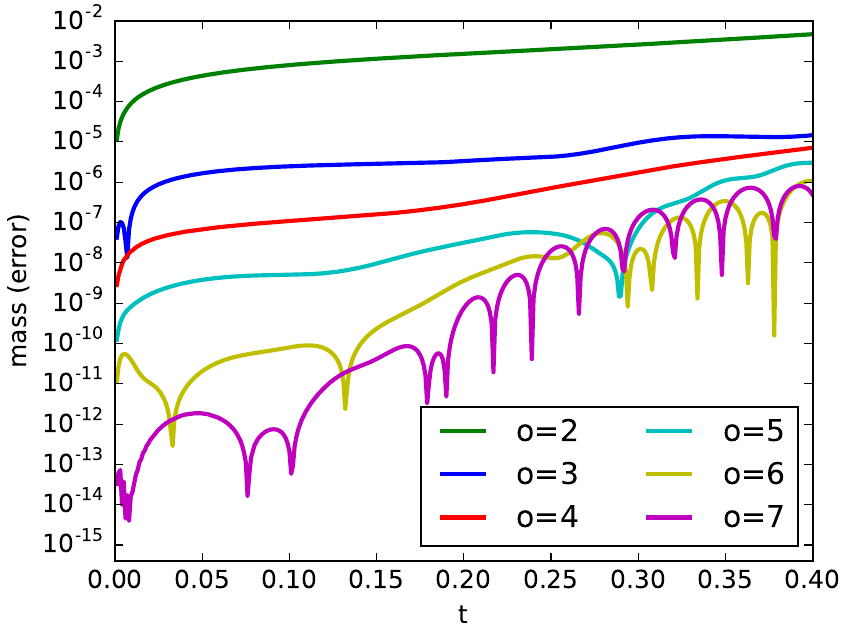}
    \caption{The numerical simulation for the KP I equation using the Schwartzian
initial value is shown at $t=0.4$ for the two-dimensional domain
(top-left) and for a one-dimensional slice with $y=0$ (top-right).
The fourth order method with a time step size of $\tau=10^{-3}$ has
been used. In addition, the error for the discontinuous Galerkin methods is shown as a function of the
    degrees of freedom \color{black}(bottom-left). The error is determined by comparing the numerical approximation to a reference solution with $o=7$ and $28672$ degrees of freedom in the $x$-direction. The error in mass is shown as a function of time (bottom-right). In the latter all configurations use $2048$ degrees of freedom in the $x$-direction. \label{fig:KPI-solution}}
\end{figure}

Second, we turn our attention to the KP I model (i.e.~$\lambda=-1$).
The numerical results for final time $t=0.4$ are shown in Figure
\ref{fig:KPI-solution}. Once again we remark that these results agree
very well with what has been reported in the literature. We find that
for low to medium precision requirements the fourth and seventh order
discontinuous Galerkin method yields the best performance, respectively.

\textcolor{black}{
    For the KP equation the mass is an invariant. However, the discontinuous Galerkin method used to solve Burgers' equation does not conserve mass exactly (see section \ref{subsec:conservation-mass}). The error in mass is shown in Figure \ref{fig:KPI-solution}, where the problem is set up such that all numerical methods use $2048$ degrees of freedom in the $x$-direction. For the fourth, fifth, sixth, and seventh order method this corresponds to an error in the infinity norm of approximately $10^{-2}$. The error in mass for these methods, however, is on the order of $10^{-6}$. Thus, even though the proposed numerical method does not conserve mass exactly, the error is very small.}

Next, let us consider the soliton solution 
\begin{equation}
u(t,x,y)=c\,\text{sech}^{2}(a(x-bt)),\label{eq:soliton}
\end{equation}
where $c=2$, $a=\epsilon\sqrt{c/2}$, and $b=2c\epsilon^{2}$. This
solution maintains the same shape in space as time goes on and we
will thus use it to test the long time behavior of our numerical method.
In the following simulation, we impose $u(0,x,y)$ as the initial
condition for the KP I equation with $\epsilon=1$. The soliton given
by equation (\ref{eq:soliton}) is only an exact solution on the entire
plane. In order to reproduce this situation as closely as possible
we will use a sufficiently large computational domain with periodic
boundary conditions.

\begin{figure}
\begin{centering}
\includegraphics[width=12cm]{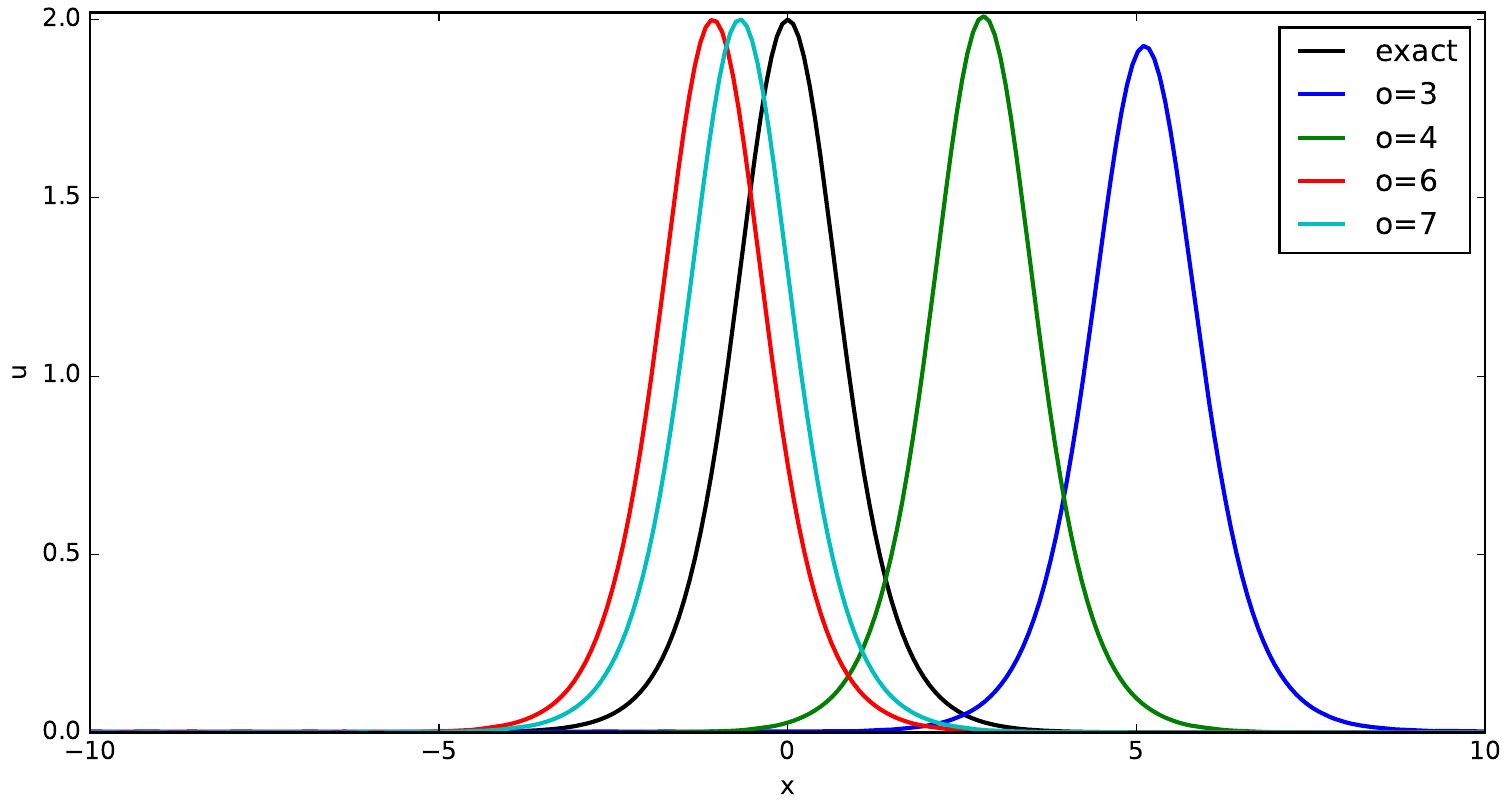}
\end{centering}
    \caption{\color{black} The slice $x\mapsto u(500,x,0)$ of the numerical approximation to equation (\ref{eq:soliton}) is shown. All numerical methods employ $512$ degrees of freedom in the $x$-direction and a time step size of $\tau=10^{-2}$ has been used. The exact solution can be computed analytically and at $t=500$ is identical to the initial value. \label{fig:soliton-solution}}
\end{figure}

For the present test we integrate the KP equation until time $t=500$.
For each of the methods we employ $512$ degrees of freedom in the
$x$-direction and a time step size $\tau=10^{-2}$ (which means that
we perform a total of $5\cdot10^{4}$ time steps). 
{\color{black} A slice of the solution (for $y=0$) at time $t=500$ is shown in Figure
\ref{fig:soliton-solution}. We remark that no spurious oscillations are visible even for the higher order methods.
}

It should be duly noted that in this case our numerical method operates far away from
the asymptotic regime (this is true for both the time integration
and the space discretization). \textcolor{black}{The described behavior is also apparent from Figure \ref{fig:soliton-solution} as the predicted speed of advection (i.e.~how far the numerical solution travels in a given time span) does not match that of the exact solution.}
As a consequence, the error in the infinity norm is large.

\textcolor{black}{Nevertheless, the proposed numerical methods (in particular, the higher order variants) still show good qualitative agreement with the exact solution. For example, the amplitude and the shape of the soliton matches that of the exact solution very well. Thus, we will consider here how well the numerical method preserves certain qualitative features of the exact solution.}

\begin{figure}
\begin{centering}
\includegraphics[width=6cm]{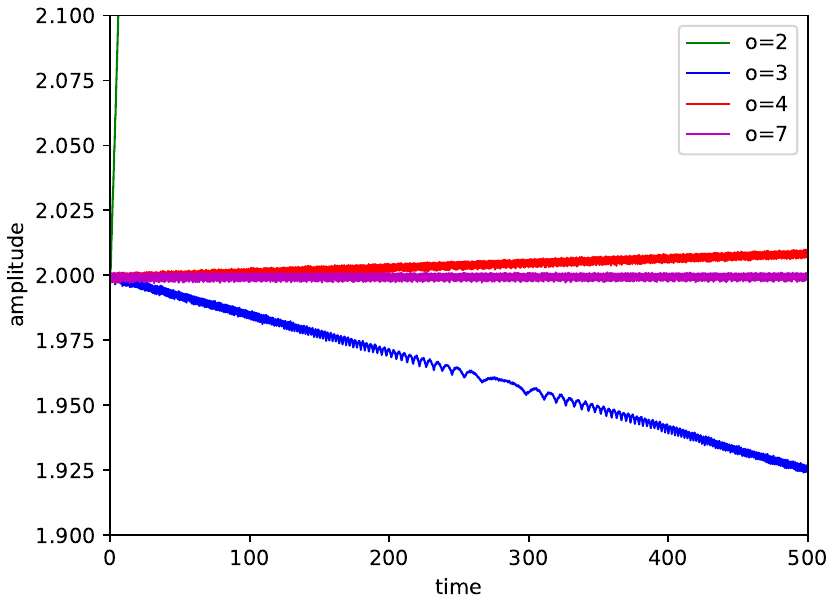}\includegraphics[width=6cm]{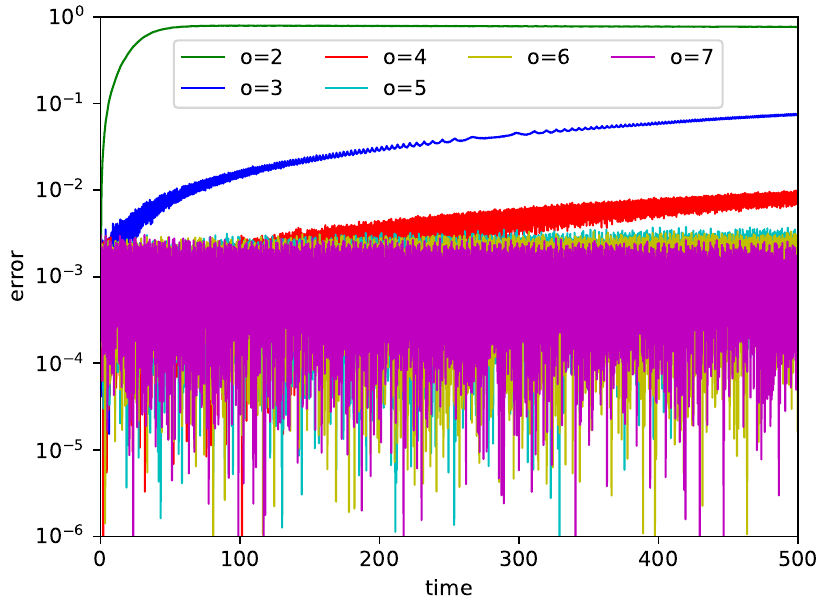}
    \begin{center}\includegraphics[width=6cm]{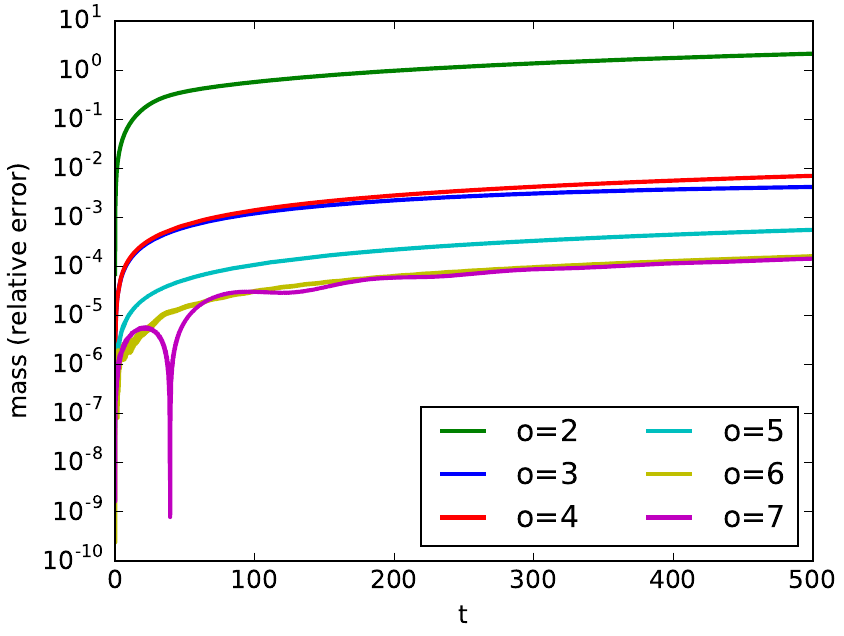}\end{center}
\par\end{centering}
    \caption{The amplitude of the numerical approximation to the soliton solution
(see equation (\ref{eq:soliton})) is shown as a function of time
    (\textcolor{black}{top-left}). In addition, the error in amplitude is plotted on a
    logarithmic scale (\textcolor{black}{top-right}). 
    \textcolor{black}{At the bottom, the relative error in mass is shown as a function of time.}
    All numerical methods employ $512$
degrees of freedom in the $x$-direction and a time step size of $\tau=10^{-2}$
has been used. \label{fig:soliton}}
\end{figure}

More specifically,
the figure of merit we are looking at here is how well the amplitude
(i.e.~the maximum) \textcolor{black}{and the mass} of the soliton is preserved. The corresponding
results are shown in Figure \ref{fig:soliton}. We observe that employing
higher order methods yields a significant decrease in the error \textcolor{black}{of both quantities}. In
particular, we see a drastic improvement going from second to third,
from third to fourth, and from fourth to fifth order, while keeping
the degrees of freedom and thus (at least approximately) the computational
cost the same.

\subsection{Performance comparison\label{subsec:comparison}}

In this section we will compare the performance of the proposed numerical
method to an approach based on exponential integrators. Exponential
integrators were evaluated for this problem in \cite{klein2011} where
they were shown to perform best compared to a range of numerical methods
(in particular, fully implicit methods and IMEX schemes). The seductive
property of these methods is that the space discretization can be
done completely with FFTs. Assuming a smooth solution this implies
superpolynomial convergence. Their disadvantage is that they treat
the nonlinearity explicitly and that (for the same number of degrees
of freedom) they are significantly more expensive compared to the
splitting approach (see \cite{einkemmer2015}).

Due to the significant difference between these two methods it is
very difficult to obtain a meaningful comparison by simply looking
at the run time of the two implementations we use in this paper. This
is in part due to the fact that in both implementations much room is
left for further optimizations (mainly from a computer science perspective).
Furthermore, we have spent significantly more time optimizing  the
splitting/discontinuous Galerkin method compared to the exponential
integrator based approach. In order to enable a meaningful comparison,
we will now develop a model that gives us a good idea of how expensive
these two algorithms are on a modern computer system.

Let us start with the splitting approach. In this case we will use
an equidistant grid in the $y$-direction in order to facilitate the
FFT that needs to be computed. However, for the $x$-direction we
will use an equidistant grid for the linear part (where we have to
take FFTs) and the discontinuous Galerkin approximation for the nonlinear
part. The former will have $n_{x}n_{y}$ degrees of freedom, while
for the latter we use $N_{x}n_{y}$ degrees of freedom (reflecting
the fact that we might need more grid points to achieve the same accuracy
for the dG method as compared to the FFT approach). We also need a
method to transfer the results of the dG computation back to the equidistant
grid (an interpolation procedure). 

In addition, the splitting (as well as the exponential integrator)
approach requires us to compute the action of the matrix exponential
$\mathrm{e}^{\tau A}u$. This is done in Fourier space and requires
the multiplication with a factor that depends on the step size $\tau$
and on the frequency. Since this factor involves an exponential, it
is rather expensive to compute. However, we can precompute these quantities
and store the result in an array. Before proceeding let us note that
this situation is, in principle, the same for the exponential integrator.
However, in the latter case more matrix functions are required, which
increases the storage cost somewhat, and care has to be taken when
evaluating the $\varphi_{k}$ functions in order to avoid large round
off errors.

Our computational model will assume that all involved operations (FFT,
interpolation, semi-Lagrangian discontinuous Galerkin scheme, array
addition and multiplication) are memory bound. That is, the performance
is dictated by how much data can be transferred to and from memory
(and not by how many arithmetic operations have to be performed).
This is a reasonable assumption on all modern computer systems. We
should note, however, that achieving this level of performance requires
an optimized implementation. In a sense this can be considered as
a (realistic) best case for both algorithms.

Using this assumption we can count the number of memory accesses per
time step that are required for the splitting/discontinuous Galerkin
scheme. The results are shown in Table \ref{tab:cost-strang}.

\begin{table}[H]
\begin{centering}
\begin{tabular}{cccc}
 & number & memory accesses & total memory accesses\tabularnewline
\hline 
FFT & $2$ & $2n_{x}n_{y}$ & $4n_{x}n_{y}$\tabularnewline
$\mathrm{e}^{\tau A}$ & $1$ & $3n_{x}n_{y}$ & $3n_{x}n_{y}$\tabularnewline
Interpolation & $1$ & $2N_{x}n_{y}$ & $2N_{x}n_{y}$\tabularnewline
dG & $1$ & $N_{x}n_{y}$ & $N_{x}n_{y}$\tabularnewline
\hline 
Strang splitting &  &  & $(7n_{x}+3N_{x})n_{y}$\tabularnewline
\end{tabular}
\par\end{centering}
\caption{The number of memory accesses required for one time step of the second
order Strang splitting/discontinuous Galerkin scheme are listed. The
number of grid points in the $y$-direction is denoted by $n_{y}$.
The number of grid points in the $x$-direction is denoted by $n_{x}$/$N_{x}$
for the equidistant/dG grid. \label{tab:cost-strang}}
\end{table}

Now, let us turn our attention to the exponential integrator of order
$2$ which can be written as
\begin{align*}
U & =\mathrm{e}^{\tau A}u^{n}+\tau\varphi_{1}(\tau A)B(u^{n})\\
u^{n+1} & =U+\tau\varphi_{2}(\tau A)(B(U)-B(u^{n}))
\end{align*}
or for our specific problem
\begin{align*}
U & =\mathrm{e}^{\tau A}u^{n}-6\tau\varphi_{1}(\tau A)(u_{x}^{n}u^{n})\\
u^{n+1} & =U-6\tau\varphi_{2}(\tau A)(UU_{x}-u_{x}^{n}u^{n}).
\end{align*}
In the following we will refer to this method as exp2. It is seductive
to just count how many FFTs, additions, etc.,~one has to compute
but this would give an overestimation of the true computational cost
(as certain operations can be combined which due to caching will reduce
the number of memory transactions). For our performance analysis we
will consider the implementation given in Table~\ref{tab:cost-exp2}. 

\begin{table}
\begin{centering}
\begin{tabular}{cccc}
step & computation & mem acc & remark\tabularnewline
\hline 
1. & $\hat{u}^{n}=\mathcal{F}(u^{n})$ & $0$ & available from previous step\tabularnewline
2. & $u_{x}^{n}=\mathcal{F}^{-1}(ik\hat{u}^{n})$ & $4n_{x}n_{y}$ & \tabularnewline
3. & $B=u_{x}^{n}\cdot u^{n}$ & $3n_{x}n_{y}$ & \tabularnewline
4. & $\hat{B}=\mathcal{F}(B)$ & $2n_{x}n_{y}$ & \tabularnewline
5. & $\hat{U}=\Phi_{0}\cdot\hat{u}^{n}+\Phi_{1}\cdot\hat{B}$ & $5n_{x}n_{y}$ & \tabularnewline
6. & $U=\mathcal{F}^{-1}(\hat{U})$ & $2n_{x}n_{y}$ & \tabularnewline
7. & $U_{x}=\mathcal{F}^{-1}(ik\hat{U})$ & $4n_{x}n_{y}$ & \tabularnewline
8. & $F=U\cdot U_{x}-B$ & $4n_{x}n_{y}$ & \tabularnewline
9. & $\hat{F}=\mathcal{F}(F)$ & $2n_{x}n_{y}$ & \tabularnewline
10. & $\hat{u}^{n+1}=\hat{U}+\Phi_{2}\cdot\hat{F}$ & $4n_{x}n_{y}$ & \tabularnewline
11. & $u^{n+1}=\mathcal{F}^{-1}(\hat{u}^{n+1})$ & $2n_{x}n_{y}$ & \tabularnewline
\hline 
Total cost &  & $32n_{x}n_{y}$ & \tabularnewline
\end{tabular}
\par\end{centering}
\caption{The different steps required for the second order exponential integrator
are listed together with the number of memory accesses (mem acc).
Note that $\mathcal{F}$ and $\mathcal{F}^{-1}$ denote the FFT and
inverse FFT, respectively. The arrays $\Phi_{0},\Phi_{1},\Phi_{2}$
contain the precomputed factors required to apply the various matrix functions
in Fourier space. \label{tab:cost-exp2}}
\end{table}

\begin{figure}[H]
\centering{}\includegraphics[width=6cm]{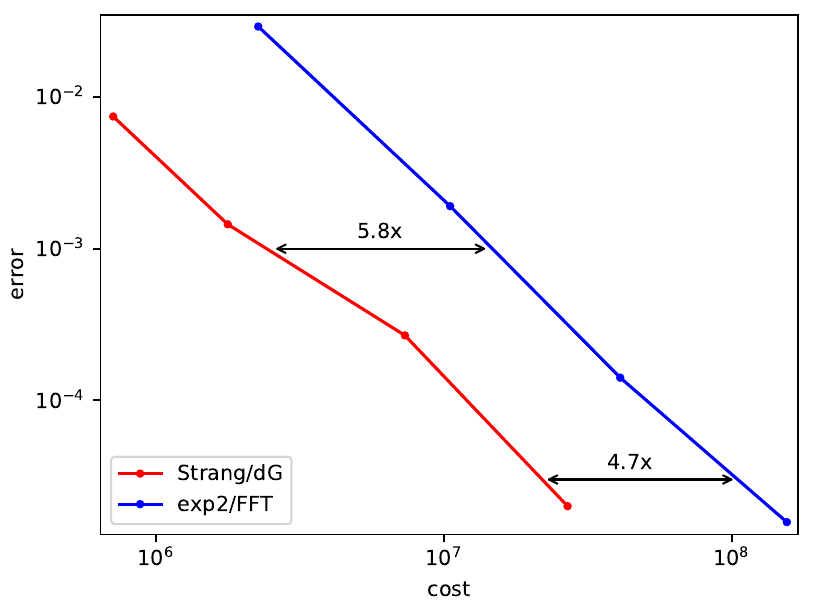}\includegraphics[width=6cm]{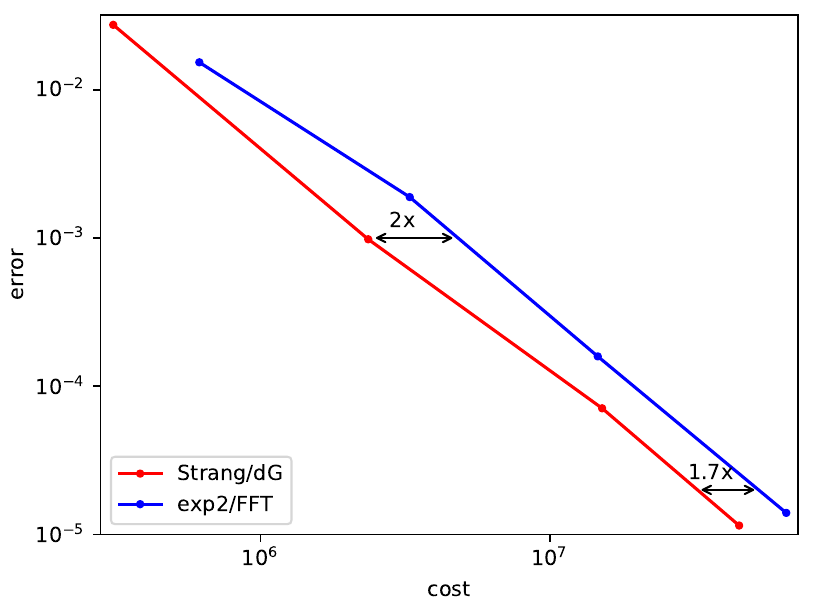}\caption{The error for the Strang splitting/discontinuous Galerkin ($o=5$)
and the exp2 scheme as a function of the computational cost is shown
for the KP I equation (left) and the KP II equation (right) using the
Schwartzian initial value (the setup from section \textcolor{black}{\ref{subsec:numerical-KP}}
is employed). The computational cost is in arbitrary units and is
determined according to the model described in this section.\label{fig:comparison}}
\end{figure}

It is clear that the Strang splitting scheme has a decisive advantage
for $N_{x}=n_{x}$ (as has already been pointed out in \cite{einkemmer2015}).
However, since the approximation error made by the Fourier transform
converges superpolynomially this is often not a very realistic assumption.
Furthermore, this does not take into account the difference in accuracy
between the splitting method and the exponential integrator. In order
to compare the two numerical methods on an equal footing we consider
for both methods a range of tolerances. At these tolerances we balance
the time and space discretization error (by choosing an appropriate
time step size and the number of degrees of freedom) and then determine
the computational cost according to the model described here. The
corresponding results for the KP equation with the Schwartzian initial
value (using precisely the setup from section \ref{subsec:numerical-KP})
are shown in Figure \ref{fig:comparison}.

We can clearly see that the Strang splitting/discontinuous Galerkin
scheme is superior for the small to medium accuracy requirements that
are often important in practice. For tighter tolerances the exp scheme
eventually becomes the method of choice (since the faster convergence
of the Fourier transform eventually overcomes the disadvantage of
the time integrator).

\section{Conclusion \& Outlook\label{sec:conclusion}}

In this paper we have proposed an arbitrary order (in space) numerical
method to solve the Kadomtsev\textendash Petviashvili equation based
on time splitting. In particular, treating the nonlinearity with a
semi-Lagrangian discontinuous Galerkin approach results in an unconditionally
stable scheme. We have demonstrated the efficiency of this numerical
method by providing a range of numerical examples. In particular,
we have compared our approach to an exponential integrator (which
has the advantage that the entire space discretization can be done
using spectral methods) and observe superior performance for small
to medium accuracy requirements (arguably the regime most important
in applications).

Although we have exclusively focused on the KP equation in this paper,
the proposed method can conceivably be extended to any partial differential
equation which combines a stiff linear part with a Burgers' type nonlinearity.
To name just a few that are of interest in the sciences: the Kawahara
equation, the (deterministic or stochastic) Kardar\textendash Parisi\textendash Zhang
equation, the viscous Burgers' equation (in particular its multi-dimensional
generalizations \cite{frisch2001}), the \mbox{Korteweg}\textendash de
Vries equation, and even the Navier\textendash Stokes equation. In
some sense the KP equation is more difficult than some of these other
models as no diffusion is included (the KP equation is a purely hyperbolic
system). 

We also consider this work as a first step towards efficiently solving
the KP equation (and other partial differential equations with higher
order derivatives) on a non-tensor-product grid. In this case the
dispersive part can not be solved using fast Fourier techniques (as
is done in this paper). Nevertheless, within the splitting approach
introduced, a preconditioned implicit method for the linear dispersion
can be combined with our semi-Lagrangian discontinuous Galerkin method
to solve Burgers' equation. We consider this as future work.

\section*{}

\bibliographystyle{plain}
\bibliography{KPdG}

\end{document}